\documentclass[letterpaper, 10 pt, conference,twocolumn]{ieeeconf}
\usepackage{amsmath,amsfonts,amssymb}
\usepackage{graphicx,epsfig,psfrag,subfigure,remark}
\usepackage{xfrac,bm}
\usepackage{pdfpages,pdfsync}
\usepackage{epstopdf}
\usepackage{pdfcolmk}
\usepackage{algorithm}
\usepackage[noend]{algpseudocode}
\usepackage[all]{xy}
\usepackage{varioref}
\usepackage{wrapfig}
\usepackage{threeparttable}
\usepackage{dcolumn}
\newcolumntype{d}{D{.}{.}{-1}}
\usepackage{nomencl}
\makeglossary
\usepackage{subfigure}
\usepackage{comment,cite}
\usepackage{xfrac}
\IEEEoverridecommandlockouts 

\newcommand{\x}{{\bm{x}}}
\newcommand{\X}{{\mathsf{X}}}

\newcommand{\f}{{\mathrm{f}}}
\newcommand{\e}{{\bm{e}}}

\newcommand{\z}{{\bm{z}}}

\newcommand{\cP}{\mathcal{P}}

\newcommand{\cC}{\mathcal{C}}

\newcommand{\cE}{\mathcal{E}}

\newcommand{\cJ}{\mathcal{J}}
\newcommand{\cK}{\mathcal{K}}
\newcommand{\cL}{\mathcal{L}}
\newcommand{\cN}{\mathcal{N}}

\newcommand{\cS}{\mathcal{S}}
\newcommand{\cU}{\mathcal{U}}

\newcommand{\mZ}{\bm{\mathbb{Z}}}

\newcommand{\mT}{T}

\newtheorem{problem}{Problem}

\renewcommand{\t}{^{\mbox{\scriptsize \mT}}}

\makeatletter
\def\BState{\State\hskip-\ALG@thistlm}
\makeatother

\renewcommand{\baselinestretch}{1.00}
\newremark{remark}{Remark}

\title{\textbf{Greedy Finite-Horizon Covariance Steering for Discrete-Time
Stochastic Nonlinear Systems Based on the Unscented Transform}}

\author{Efstathios Bakolas \and Alexandros Tsolovikos \thanks{E. Bakolas
is an Associate Professor and A. Tsolovikos is a PhD student in the Department of Aerospace Engineering
and Engineering Mechanics, The University of Texas at Austin,
Austin, Texas 78712-1221, USA, Email: bakolas@austin.utexas.edu; tsolovikos@austin.utexas.edu.
This work was supported in part by the National Science Foundation
(award no. CMMI-1753687).} }

\begin{document}

\maketitle

\begin{abstract}
In this work, we consider the problem of steering the first two
moments of the uncertain state of a discrete-time nonlinear
stochastic system to prescribed goal quantities at a given final
time. In principle, the latter problem can be formulated as a
density tracking problem, which seeks for a feedback policy that
will keep the probability density function of the state of the
system close, in terms of an appropriate metric, to the goal
density. The solution to the latter infinite-dimensional problem can
be, however, a complex and computationally expensive task. Instead,
we propose a more tractable and intuitive approach which relies on a
greedy control policy. The latter control policy is comprised of the
first elements of the control policies that solve a sequence of
corresponding linearized covariance steering problems. Each of these
covariance steering problems relies only on information available
about the state mean and state covariance at the current stage and
can be formulated as a tractable (finite-dimensional) convex
program. At each stage, the information on the state statistics is
updated by computing approximations of the predicted state mean and
covariance of the resulting closed-loop nonlinear system at the next
stage by utilizing the (scaled) unscented transform. Numerical
simulations that illustrate the key ideas of our approach are also
presented.
\end{abstract}

\section{Introduction}\label{s:intro}
This paper deals with the finite-horizon covariance steering problem
for discrete-time stochastic nonlinear (DTSN) systems. In
particular, we consider the problem of steering the first moment
(mean) and the second central moment (covariance) of the uncertain
state of a DTSN system to desired quantities at a given (finite)
terminal time. We will refer to the latter problem as the nonlinear
covariance steering problem to emphasize the fact that it is the
steering of the state covariance that constitutes the most
challenging and less studied part of this stochastic control problem
(steering the state mean essentially corresponds to a standard, but
not necessarily trivial, controllability problem). Perhaps, one of
the most natural approaches to address nonlinear covariance steering
problems would be to place them under the umbrella of PDE tracking
problems in which one tries to minimize the distance of the
probability density of the state, which evolves in space and time in
accordance with the Fokker Planck partial differential equation
(PDE), from a desired terminal density
function~\cite{p:FLEIG201913}. The solution to the latter
infinite-dimensional optimization problem, however, can be a very
complex task in general. In this work, we will employ a more
practical approach that relies on the solution of a sequence of
linearized steering problems which are in turn reduced to tractable
convex optimization problems.

\textit{Literature Review:} In the special case of linear Gaussian
systems, that is, stochastic linear systems subject to Gaussian
white noise, the covariance steering problem corresponds to a
distribution steering problem, in the sense that the mean and
covariance of the terminal state uniquely determine the (Gaussian)
probability distribution of the latter state. Infinite-horizon
covariance steering (also known as covariance control) problems for
both continuous-time and discrete-time Gaussian systems have been
studied extensively by Skelton and his co-authors in a series of
papers (see, for instance,
\cite{p:skeltonIJC,p:skeltonTAC,p:yasuda1993,p:Grig97,p:levy1997discrete}).
The finite-horizon problem for the continuous time case was recently
revisited and studied in detail
in~\cite{p:georgiou15A,p:georgiou15B} whereas the same problem for
the discrete-time case was studied
in~\cite{p:bakCDC16,p:PT2017,p:BAKOLAS2018}. The previous references
assume perfect state information (that is, the realization of the
state process at each state can be measured perfectly). Covariance
control problems in the case of incomplete and imperfect state
information have been studied in
\cite{p:EBACC17,p:bakCDC2018,p:bakTAC2019}. 
Nonlinear density steering problems for feedback linearizable
nonlinear systems were recently studied in~\cite{p:caluya2019}. An
iterative covariance steering algorithm for nonlinear systems based
on a simple linearization of the system dynamics along reference
state and input trajectories can be found in~\cite{p:ridderhof2019}. Stochastic nonlinear model predictive control with probabilistic
constraints can be found in~\cite{p:mesbah2014,p:liu2014,p:SEHR2017,p:PAULSON2018}.

\textit{Main Contribution:} In this work, we propose a greedy, yet
practical and intuitive, solution approach to the nonlinear
covariance steering problem. The proposed approach consists of three
key steps which are applied iteratively. In the first step, we
linearize the system dynamics around the current state of the system
(rather than along a reference trajectory). The particular
linearization scheme relies on information available at the current
stage and in particular, knowledge of approximations of the mean and
covariance of the current state of the system. At each new stage, a
new linearization will be computed to account for the new
information that becomes available at that stage. We refer to the
first step as the \textit{recursive linearization} step (RL step).
In the second step, we compute a feedback control policy (sequence
of feedback control laws) that solves a relevant linear, Gaussian
covariance steering problem based on available approximations of the
current state mean and covariance and the linear state space model
computed at the LN step. The latter policy can be computed in
real-time by means of tractable convex optimization techniques by
leveraging the results of our previous work in covariance steering
problems for Gaussian linear
systems~\cite{p:bakCDC16,p:PT2017,p:BAKOLAS2018}. From the computed
policy, only the first control law is executed at each stage. We
refer to the latter step as the \textit{linearized Gaussian
covariance steering} step (LGCS step).

In the third step, we compute approximations of the
one-stage-predictions of the state mean and covariance of the
closed-loop system that results by applying the feedback control
policy computed at the LGCS step. To compute these approximations,
we employ the (scaled) unscented
transform~\cite{p:julier2002,p:julier2004}. The latter transform
relies on the propagation of a small number of points, which are
known as ``sigma points,'' in future stages. These points are
selected in a deterministic way such that their mean and variance
are compatible with prior information~\cite{p:julier2002}. The
predicted state mean and covariance of the closed-loop system
determine a Gaussian (or normal) approximation of the (predicted)
state statistics of the next state. For this reason, we shall refer
to the latter step as the \textit{predictive normalization} step (PN
step). This three-step process is repeated iteratively until the
final stage, when it is expected that the (terminal) state mean and
covariance are sufficiently close to the goal quantities.

The previously described iterative process corresponds to an on-line
(or real-time) greedy control policy for nonlinear covariance
steering. Because predictions of the state statistics in this
approach do not go beyond the next stage, there cannot be explicit
performance considerations as in a typical model predictive control
approach~\cite{p:grune2017}. Instead, the emphasis of the proposed
greedy approach is placed on satisfying as closely as possible the
boundary conditions (by steering the state mean and state covariance
to desired prescribed quantities).

\textit{Structure of the paper:} The rest of the paper is organized
as follows. In Section~\ref{s:form}, we formulate the nonlinear
covariance steering problem. A greedy algorithm for the solution to
the latter problem is presented in Section~\ref{s:greedy}.
Furthermore, we present numerical simulations in
Section~\ref{s:simu} and we conclude the paper with a number of
remarks and directions for future research in Section~\ref{s:concl}.

\section{Problem Formulation}\label{s:form}

\subsection{Notation}
We denote by $\mathbb{R}^n$ the set of $n$-dimensional real vectors.
Given integers $\alpha, \beta$ with $\alpha \leq \beta$, we denote
by $[\alpha, \beta]_d$, the discrete interval from $\alpha$ to
$\beta$. We denote by $\mathbb{E}[\cdot]$ the expectation operator.
Given a random vector $x$, we denote by $\mathbb{E}[x]$ its mean and
by $\mathrm{Cov}[x]$ it covariance, where $\mathrm{Cov}[x]:=
\mathbb{E}[ (x- \mathbb{E}[x])(x- \mathbb{E}[x])\t ]$. The space of
real symmetric $n \times n$ matrices will be denoted by
$\mathbb{S}_n$. Furthermore, we will denote the convex cone of
$n\times n$ (symmetric) positive semi-definite and (symmetric)
positive definite matrices by $\mathbb{S}^{+}_n$ and
$\mathbb{S}^{++}_n$, respectively. Finally, we write
$\mathrm{bdiag}(A_1,$ $\dots, A_\ell)$ to denote the block diagonal
matrix formed by the matrices $A_i$, $i\in \{1,\dots, \ell\}$.

\subsection{Problem setup}
We consider the following discrete-time nonlinear stochastic system
\begin{equation}\label{eq:nonlinear}
x(t+1) = f(x(t), u(t)) + w(t),
\end{equation}
for $t \in [0,N-1]_{d}$, where $N$ is a positive integer, and $x(0)
=x_0$, where $x_0$ is a random vector with $\mathbb{E}[x_0] = \mu_0$
and $\mathrm{Cov}[x_0] = \Sigma_0$, with $\mu_0 \in \mathbb{R}^n$
and $\Sigma_0 \in \mathbb{S}^{++}_n$ be given quantities.
Furthermore, $f(\cdot)$ is a $C^1$ function. In addition, $x_{0:N}
:= \{ x(t) \in \mathbb{R}^n: t\in [0,N]_{d} \}$ and $u_{0:N-1} := \{
u(t) \in \mathbb{R}^m:~t\in [0,N-1]_{d}\}$ correspond to the state
and input (random) processes, respectively. In addition, $w_{0:N-1}
:= \{ w(t) \in \mathbb{R}^n:~t\in [0,N-1]_{d}\}$ corresponds to the
noise process which is assumed to be a sequence of independent and
identically distributed random variables with
\begin{align}\label{eq:w1}
\mathbb{E}\left[ w(t) \right] & = 0, &\mathbb{E}\left[ w(t)
w(\tau)\t \right] & = \delta(t,\tau) W_t,
\end{align}
for all $t,\tau\in [0,N-1]_{d}$, where $W_t\in\mathbb{S}_n^{+}$,
and $\delta(t,\tau):=1$, when $t=\tau$, and $\delta(t,\tau):=0$,
otherwise. Furthermore, $x_0$ is independent of $w_{0:N-1}$, that
is,
\begin{subequations}
\begin{align}\label{eq:x0w}
\mathbb{E}\left[ x_0  w(t)\t  \right] & = 0,~~~~\mathbb{E}\left[
w(t) x_0\t  \right] = 0,
\end{align}
\end{subequations}
for all $t \in [0,N-1]_d$. Finally, throughout this paper we assume
that we have perfect state information, that is, at each stage $t$,
the realization $x(t)$ of the state process is perfectly known
(measured).

Because the system given in \eqref{eq:nonlinear} is nonlinear, even
if the initial state is drawn from a normal distribution and the
noise is white Gaussian, it is not guaranteed that the state at
future stages will remain Gaussian. For this reason, it is not
meaningful to require that the terminal state of the system should
be steered to a prescribed normal distribution as in the standard
formulation of a finite-horizon covariance steering problem for
Gaussian linear systems. A more practical approach would be to
require that the state mean and covariance of the nonlinear system
attain (exactly or approximately) prescribed quantities. In
particular, let us denote by $\mu_x(t)$ and $\Sigma_x(t)$ the state
mean and covariance at stage $t$, that is,
\begin{align}
\mu_x(t) := \mathbb{E}[x(t)], ~~~~ \Sigma_x(t) & : =
\mathrm{Cov}[x(t)]. 
\end{align}
The class of admissible control policies is taken to be the set of
sequences of control laws that are measurable functions of the
realization of the current state of the system. Then, the nonlinear
covariance steering problem can be formulated as follows:
\begin{problem}\label{problem1}
Let $\mu_0, \mu_\f\in\mathbb{R}^n$ and $\Sigma_0, \Sigma_\f \in
\mathbb{S}^{++}_n$ be given. Find a control policy $\pi : = \{
\kappa(t,\cdot): t\in [0,N-1]_d) \}$ that will steer the system
\eqref{eq:nonlinear} from $x(0) = x_0$ with $\mathbb{E}[x_0]=\mu_0$
and $\mathrm{Cov}[x_0] = \Sigma_0$ to a terminal state $x(N)$ with
\begin{equation}
\mu_x(N) =\mu_\f,~~~~ (\Sigma_\f - \Sigma_x(N)) \in
\mathbb{S}^{+}_n.
\end{equation}
\end{problem}

\subsection{Collection of Finite-Horizon Linearized Covariance Steering Problems}
Next, we associate the DTSN system~\eqref{eq:nonlinear} at stage
$t=k \in [0,N-1]_d$ with a discrete-time stochastic linear system.
The latter system corresponds to a linearization of the DTSN system
around a given point $(\mu_k, \nu_k)\in \mathbb{R}^n \times
\mathbb{R}^m$ which is given by
\begin{align}\label{eq:linear}
z(t+1) & = A_k (z(t) - \mu_k) + B_k (u(t) - \nu_{k}) \nonumber \\
&~~~~~~ + r_k + w(t),
\end{align}
for $t\in[k, N-1]_d$ and $z(k) = \mathsf{z}_k$, with
$\mathbb{E}[\mathsf{z}_k] = \mu_k$ and $\mathrm{Cov}[\mathsf{z}_k] =
\Sigma_k$, where $\mu_k \in \mathbb{R}^n$, $\Sigma_k \in
\mathbb{S}^{++}_n$, and $\nu_{k} \in \mathbb{R}^m$. 
In addition, it is assumed that
\begin{subequations}
\begin{align}\label{eq:z0w}
\mathbb{E}\left[ \mathsf{z}_k  w(t)\t  \right] & =
0,~~\mathbb{E}\left[ w(t)\mathsf{z}_k\t  \right] = 0,
\end{align}
\end{subequations}
for all $t \in [k,N-1]_d$. Furthermore, $A_k$ and $B_k$ are constant
(time-invariant) matrices whereas $r_k$ is a constant vector, and in
particular,
\begin{subequations}
\begin{align}\label{eq:AkBk}
A_k & := \frac{\partial }{\partial x} f(x, u) \Big|_{ \substack{x = \mu_k \\
u = \nu_k}},~~B_k := \frac{\partial }{\partial u} f(x, u) \Big|_{
\substack{x = \mu_k\\ u= \nu_k }}, \\
r_k & := f(\mu_k, \nu_k). \label{eq:rk}
\end{align}
\end{subequations}
We can equivalently write \eqref{eq:linear} in a slightly more
compact form as follows:
\begin{align}\label{eq:linearB}
z(t+1) & = A_k z(t) + B_k u(t) + d_k + w(t),
\end{align}
where $d_k := -A_k \mu_k - B_k \nu_{k} + r_k$. 

We will refer to the latter linear model as the $k$-th linearized
state space model ($t=k$ corresponds to the initial stage). It is
worth noting that the triple $(A_k, B_k, r_k)$ remains constant
throughout the whole horizon $[k,N-1]_d$. However, for a different
$k$, one obtains a different linearized system with a different but
time-invariant triplet $(A_k, B_k, r_k)$. Therefore,
\eqref{eq:linear} describes essentially a collection of $N-k$
different (one for each $k$) time-invariant systems. An implicit
assumption here is that the pair $(A_k, B_k)$ is controllable.

An alternative linear model can be derived if one linearizes the
DTSN system \eqref{eq:nonlinear} around $(\mu_\f,\nu_\f) \in
\mathbb{R}^n \times \mathbb{R}^m$, where $\mu_\f$ is the goal mean
of the terminal state and $\nu_\f$ is such that $\mu_\f = f(\mu_\f,
\nu_\f)$. Then, $A_k = A_\f$, $B_k = B_\f$, and $r_k = r_\f$, where
\begin{subequations}
\begin{align}
A & := \frac{\partial }{\partial x} f(x, u) \Big|_{ \substack{x = \mu_\f \\
u = \nu_\f}},~~B := \frac{\partial }{\partial u} f(x, u) \Big|_{
\substack{x = \mu_\f \\ u= \nu_\f }} \\
r_k & = f(\mu_\f, \nu_\f),
\end{align}
\end{subequations}
for all $k \in [0,N-1]$. The previous linearization assumes that the
DTSN system operates ``near'' the terminal target point $(\mu_\f,
\nu_\f) \in \mathbb{R}^n \times \mathbb{R}^m$. 

Note that both of the previously described linearized models are
different from the one obtained after linearizing a nonlinear system
around a given pair of reference state and input sequences
$\bar{z}_{0:N} := \{\bar{z}(t): t\in [0,N]_d\}$ and $\bar{u}_{0:N-1}
:= \{\bar{u}(t): t\in [0,N-1]_d\}$, respectively, as is proposed,
for instance, in \cite{p:ridderhof2019}. In the latter case, one
would consider a single time-varying linearized system described by
the following equation:
\begin{equation}\label{eq:linear0}
z(t+1) = A(t) z(t) + B(t) u(t) + r(t) + w(t),
\end{equation}
for $t\in[0, N-1]_d$, where $A(t)$, $B(t)$, and $r(t)$ are
time-varying matrices which are defined as follows:
\begin{align*}
A(t) & := \frac{\partial }{\partial x}f(x, u) \Big|_{
\substack{x=\bar{z}(t) \\  u=\bar{u}(t)}}, ~~~ B(t) :=
\frac{\partial }{\partial u}f(x, u) \Big|_{\substack{x=\bar{z}(t) \\
u = \bar{u}(t)}},\\
r(t) & : = f(\bar{z}(t), \bar{u}(t)) - A(t) \bar{z}(t) - B(t)
\bar{u}(t),
\end{align*}
for all $t\in[0,N-1]_d$.

However, finding a reference state sequence $\bar{z}_{0:N}$ and a
corresponding (compatible) reference input sequence
$\bar{u}_{0:N-1}$ may be a non-trivial task.
In particular, the reference state sequence should
satisfy the desired boundary conditions whereas the reference input
sequence should generate the corresponding reference state sequence.

Next, we will formulate a linearized covariance steering problem for
the system described in \eqref{eq:linear} for a given $k\in
[0,N-1]_d$. The class $\cU$ of admissible control policies for the
latter problem will consist of sequence of control laws $\{\phi_k(t,
\cdot ):~ t\in[k,N-1]_d \}$, where
\begin{equation}\label{eq:phi}
\phi_k(t, z) = \upsilon_k(t) + K_k(t) z,~~~~t\in[k,N-1]_d.
\end{equation}

We next formulate the linearized covariance steering based on
information available at stage $t=k$.
\begin{problem}[$k$-th linearized covariance steering problem] \label{problem2}
Let $\mu_k, \mu_\f\in\mathbb{R}^n$ and $\Sigma_k, \Sigma_\f \in
\mathbb{S}^{++}_n$ be given. Among all admissible control policies
$\varpi_k : = \{ \phi_k(k,\cdot), \dots, \phi_k(N-1,\cdot) \} \in
\cU$, where $\phi_k(t, \cdot)$ satisfies \eqref{eq:phi} for $t\in
[k,N-1]_d$, find a control policy $\varpi_k^{\star}$ that minimizes
the following performance index
\begin{align}\label{eq:cost}
J_k(\varpi_k) & := 
\mathbb{E}\big[\sum_{t = k}^{N-1} \phi_k(t,z(t))\t \phi_k(t,z(t))
\big]
\end{align}
subject to the recursive dynamic constraints \eqref{eq:linear} and
the following boundary conditions:
\begin{subequations}
\begin{align}
\mathbb{E}[\mathsf{z}_k] & = \mu_k, & \mathrm{Cov}[ \mathsf{z}_k ] = \Sigma_k,\\
\mathbb{E}[z(N)] & =\mu_\f, & (\Sigma_\f - \mathrm{Cov}[z(N)] ) \in
\mathbb{S}^{+}_n.
\end{align}
\end{subequations}

\end{problem}
\begin{remark}
The choice of the performance index is to ensure that the control
input will have finite energy and thus avoid excessive actuation (as
we have already mentioned, performance considerations are not of
primary interest in this work). Problem~\ref{problem2} does not
correspond to a standard finite-horizon linear quadratic Gaussian
(LQG) problem due to the presence of the (non-standard) terminal
positive semi-definite constraint $(\Sigma_\f - \mathrm{Cov}[
\Sigma_z(N) ]) \in \mathbb{S}^{+}_n$. Although we do not explicitly
consider state or input constraints in the formulation of
Problem~\ref{problem2}, we will present an optimization-based
solution which is also applicable to more general problem
formulations including those with input and/or state constraints as
in \cite{p:BAKOLAS2018}.
\end{remark}

\begin{remark}
Note that finding a policy $\varpi_k$ that solves
Problem~\ref{problem2} is equivalent to finding a sequence $\{
(\upsilon_k(t),K_k(t)):~t\in[k,N-1]_d \}$. The main idea of the
proposed solution approach is that the first control law of the
control policy that solves the $k$-th linearized covariance steering
problem (Problem~\ref{problem2}) can be used as the control law
corresponding to the stage $t=k$ from the control policy that is a
candidate solution to the nonlinear covariance steering problem
(Problem~\ref{problem1}). Later on, we will see that this idea will
have to be applied iteratively in the sense that the data of the
linear covariance problem will change at each stage, and
consequently, the corresponding feedback policy has to be updated
accordingly to reflect the new information available. 
\end{remark}

\subsection{Solution to the $k$-th Linearized Covariance Steering}
Next, we will present the main steps of the solution to the $k$-th
linearized covariance steering problem (Problem~\ref{problem2}). To
this aim, Eq.~\eqref{eq:linear} can be written in compact form
as follows:
\begin{align}\label{eq:bmz}
\bm{z} = \mathbf{G}^k_z \mathsf{z}_k + \mathbf{G}^k_{u} \bm{u} +
\mathbf{G}^k_{w}( \bm{w} + \bm{d}_k),
\end{align}
where
\begin{align*}
\bm{z} & := [z(k)\t, \dots , z(N)\t]\t, ~~
\bm{u} := [u(k)\t, \dots , u(N-1)\t]\t \\
\bm{w} & := [w(k)\t, \dots , w(N-1)\t]\t ~~
\bm{d}_k := [d_k\t, \dots , d_k\t]\t.
\end{align*}
In addition, $\mathbf{G}^k_{u}$, $\mathbf{G}^k_{w}$, and
$\mathbf{G}^k_{z}$ are defined as follows:
\begin{align*}
& \mathbf{G}^k_{u} := \begin{bmatrix} 0 & 0 & \dots & 0 \\
B_k &  0 & \dots & 0 \\
A_k B_k & B_k & \dots & 0 \\
\vdots & \vdots & \dots & \vdots \\
A_k^{N-1-k} B_k & A_k^{N-2-k} B_k &
\dots & B_k\end{bmatrix},\\
&\mathbf{G}^k_{w} :=\begin{bmatrix} 0 & 0 & \dots & 0 \\
I & 0 & \dots & 0 \\
A_k & I & \dots & 0 \\
\vdots & \vdots & \dots & \vdots \\
A_k^{N-1-k} & A_k^{N-2-k} & \dots & I
\end{bmatrix}, \\
& \mathbf{G}^k_{z} :=
\begin{bmatrix}
I & A_k\t & \dots & \big(A_k^{N-k}\big)\t
\end{bmatrix}\t.
\end{align*}

In view of \eqref{eq:phi}, an admissible control sequence can be
written compactly a follows:
\begin{equation}\label{eq:bmucloop0}
\bm{u} = \bm{\cK}_{k} \bm{z} + \bm{\upsilon}_k,
\end{equation}
where 
\begin{align*}
\bm{\cK}_{k} & := 
[\mathrm{bdiag}( K_{k}(k), \dots, K_{k}(N-1)),~0], \\
\bm{\upsilon}_{k} & := [ \upsilon_{k}(k)\t, \dots,
\upsilon_{k}(N-1)\t]\t.
\end{align*}

Consequently, after plugging \eqref{eq:bmucloop0} into
\eqref{eq:bmz}, we can express the closed-loop dynamics in compact
form as follows:
\begin{align}\label{eq:bmzcloop}
\bm{z} = \mathbf{T}_z^k \mathsf{z}_k + \mathbf{T}_{\upsilon}^k
\bm{\upsilon}_k + \mathbf{T}_{w}^k (\bm{w} + \bm{d}_k)
\end{align}
where
\begin{subequations}
\begin{align}\label{eq:TFzcloop}
\mathbf{T}_{z}^k & := (I - \mathbf{G}^k_{u}
\bm{\cK}_{k})^{-1}  \mathbf{G}^k_z \\
\mathbf{T}_{\upsilon}^k & := (I - \mathbf{G}^k_{u}
\bm{\cK}_{k})^{-1}  \mathbf{G}^k_{u} \label{eq:TFzcloop2} \\
\mathbf{T}_{w}^k & := (I - \mathbf{G}^k_{u} \bm{\cK}_{k})^{-1}
\mathbf{G}^k_{w}. \label{eq:TFzcloop3}
\end{align}
\end{subequations}
Note that the matrix $(I - \mathbf{G}^k_{u} \bm{\cK}_{k})$
corresponds to a block lower triangular matrix whose diagonal blocks
are equal to the identity matrix~(for more details, the reader may
refer to~\cite{p:skaf2010,p:BAKOLAS2018}). Thus, $(I -
\mathbf{G}^k_{u} \bm{\cK}_{k})^{-1}$ is well-defined.

In view of equation \eqref{eq:bmzcloop}, \eqref{eq:bmucloop0}
becomes
\begin{align}\label{eq:bmucloop2}
\bm{u} = \mathbf{H}_z^k \mathsf{z}_k + \mathbf{H}_{\upsilon}^k
\bm{\upsilon}_k + \mathbf{H}_{w}^k (\bm{w} + \bm{d}_k),
\end{align}
where
\begin{subequations}
\begin{align}\label{eq:HFucloop}
\mathbf{H}_{z}^k & := \bm{\cK}_{k}(I - \mathbf{G}^k_{u}
\bm{\cK}_{k})^{-1}  \mathbf{G}^k_z \\
\mathbf{H}_{\upsilon}^k & := I + \bm{\cK}_{k}(I - \mathbf{G}^k_{u}
\bm{\cK}_{k})^{-1}  \mathbf{G}^k_{u} \label{eq:HFucloop2}\\
\mathbf{H}_{w}^k & := \bm{\cK}_{k}(I - \mathbf{G}^k_{u}
\bm{\cK}_{k})^{-1}  \mathbf{G}^k_{w}. \label{eq:HFucloop3}
\end{align}
\end{subequations}

After plugging \eqref{eq:bmucloop2} in \eqref{eq:cost}, one can
obtain an expression for the (predicted) cost as a function of the
decision variables $\bm{\cK}_{k}$ and $\bm{\upsilon}_{k}$. In
particular,
\begin{align}\label{eq:costKu}
J_k(\varpi_k)
& = \mathbb{E}\big[ \bm{u}\t \bm{u} \big] = \mathrm{trace}\big(
\mathbb{E}\big[ \bm{u} \bm{u}\t \big] \big) \nonumber \\
%
& = \mathrm{trace}\big( \mathbb{E}\big[ \big( \mathbf{H}_z^k
\mathsf{z}_k + \mathbf{H}_{\upsilon}^k \bm{\upsilon}_k +
\mathbf{H}_{w}^k ( \bm{w} + \bm{d}_k) \big)\nonumber \\
&~~\qquad~~~ \times \big( \mathbf{H}_z^k \mathsf{z}_k +
\mathbf{H}_{\upsilon}^k \bm{\upsilon}_k + \mathbf{H}_{w}^k ( \bm{w}
+ \bm{d}_k )
\big)\t \big] \big) \nonumber \\
& =: \tilde{J}_k(\bm{\cK}_{k},\bm{\upsilon}_k).
\end{align}
It follows readily that
\begin{align}
\tilde{J}_k(\bm{\cK}_{k}, \bm{\upsilon}_k) & = \mathrm{trace}\big(
\mathbf{H}_z^k (\Sigma_k + \mu_k\mu_k\t) (\mathbf{H}_z^k)\t \nonumber \\
&~~~~~~~ + 2 \mathbf{H}_z^k \mu_k \bm{\upsilon}_k\t
(\mathbf{H}_{\upsilon}^k)\t  + 2 \mathbf{H}_z^k \mu_k \bm{d}_k\t
(\mathbf{H}_{w}^k)\t \nonumber \\
&~~~~~~~ + \mathbf{H}_{\upsilon}^k \bm{\upsilon}_k \bm{\upsilon}_k\t
(\mathbf{H}_{\upsilon}^k)\t + 2 \mathbf{H}_{\upsilon}^k
\bm{\upsilon}_k \bm{d}_k\t
(\mathbf{H}_{w}^k)\t \nonumber \\
&~~~~~~~ + \mathbf{H}_{w}^k (\mathbf{W}_{k:N-1} +
\bm{d}_k\bm{d}_k\t) (\mathbf{H}_{w}^k)\t \big),
\end{align}
where $\mathbf{W}_{k:N-1}:= \mathrm{bdiag}(W_k,\dots, W_{N-1})$. In
the previous derivation, we have used the available information
about the statistics of $\mathsf{z}_k$ and in particular, that
$\mathbb{E}[\mathsf{z}_k] = \mu_k$, $\mathbb{E}[\mathsf{z}_k
\mathsf{z}_k\t] = \Sigma_k + \mu_k \mu_k\t$. 

Next, we express the terminal constraints in terms of the decision
variables $(\bm{\cK}_k, \bm{\upsilon}_k)$. In particular, we have
\begin{align}\label{eq:meancloop}
\mathbb{E}[z(N)] & =  \mathbb{E}[\mathbf{P}_N \bm{z}] = \mathbf{P}_N \mathbb{E}[\bm{z}] \nonumber \\
& = \mathbf{P}_N \big( \mathbf{T}^k_{z} \mu_k +
\mathbf{T}^k_{\upsilon} \bm{\upsilon}_k + \mathbf{T}^k_{w}
\bm{d}_k \big) \nonumber \\
& =: \mathfrak{f}(\bm{\cK}_k, \bm{\upsilon}_k),
\end{align}
where $\mathbf{P}_N := [0, \dots, 0,~I]$. Therefore, the constraint
$\mathbb{E}[z(N)] =\mu_\f$ can be written as follows:
\begin{equation}\label{eq:C1}
C_1(\bm{\cK}_k, \bm{\upsilon}_k) =0,~~~C_1(\bm{\cK}_k,
\bm{\upsilon}_k) := \mathfrak{f}(\bm{\cK}_k, \bm{\upsilon}_k) -
\mu_\f,
\end{equation}
where $\mathfrak{f}(\bm{\cK}_k, \bm{\upsilon}_k)$ is given in
\eqref{eq:meancloop}. Furthermore, we have that
\begin{align}\label{eq:covarbnd}
\mathrm{Cov}[ z(N) ] = \mathbb{E}[ z(N) z(N)\t ] - \mu_\f \mu_\f\t,
\end{align}
where
\begin{align}\label{eq:EzztN}
\mathbb{E}[ z(N) z(N)\t ] & = \mathbf{P}_N \mathbb{E}\big[ \big( \mathbf{T}_z^k 
\mathsf{z}_k + \mathbf{T}_{\upsilon}^k \bm{\upsilon}_k +
\mathbf{T}_{w}^k ( \bm{w} + \bm{d}_k ) \big)\nonumber \\
& ~~~ \times \big( \mathbf{T}_z^k \mathsf{z}_k +
\mathbf{T}_{\upsilon}^k \bm{\upsilon}_k + \mathbf{T}_{w}^k ( \bm{w}
+ \bm{d}_k ) \big)\t \big] \mathbf{P}_N\t \nonumber \\
& =: \mathfrak{g}(\bm{\cK}_k, \bm{\upsilon}_k).
\end{align}
Therefore, the terminal state covariance constraint: $(\Sigma_\f -
\Sigma_z(N)) \in \mathbb{S}_n^{+}$, can be written as the following
positive semi-definite constraint:
\begin{subequations}
\begin{align}\label{eq:C2}
C_2(\bm{\cK}_{k}, \bm{\upsilon}_k) & \in \mathbb{S}_n^{+},\\
C_2(\bm{\cK}_{k},\bm{\upsilon}_k) & := \Sigma_\f -
\mathfrak{g}(\bm{\cK}_k, \bm{\upsilon}_k) + \mu_\f \mu_\f\t,
\label{eq:C2b}
\end{align}
\end{subequations}
where $\mathfrak{g}(\bm{\cK}_k, \bm{\upsilon}_k)$ is defined in
\eqref{eq:EzztN}.

\begin{problem}\label{problem:optiz}
Find a pair $ (\bm{\cK}_{k}^{\star}, \bm{\upsilon}_k^{\star})$ that
minimizes the predicted cost $\tilde{J}_k(\bm{\cK}_{k},
\bm{\upsilon}_k)$ subject to the constraints:
\begin{align}\label{eq:C1C2}
C_1(\bm{\cK}_k, \bm{\upsilon}_k) = 0,~~~~ C_2(\bm{\cK}_{k},
\bm{\upsilon}_k) \in \mathbb{S}^{+}_n,
\end{align}
where $C_1(\bm{\cK}_k, \bm{\upsilon}_k)$ and $C_2(\bm{\cK}_{k},
\bm{\upsilon}_k)$ are defined in \eqref{eq:C1} and \eqref{eq:C2b},
respectively.
\end{problem}

Problem~\ref{problem:optiz} is not convex as is explained in
\cite{p:BAKOLAS2018}. One can associate it, however, with a convex
program by applying suitable transformations to the pair of decision
variables $( \bm{\cK}_k, \bm{\upsilon}_k)$ in order to obtain a new
pair of decision variables, $(\bm{\cL}_k, \bm{\nu}_k)$, which are
defined as follows~\cite{p:skaf2010}:
\begin{subequations}
\begin{align}
\bm{\cL}_k & := \bm{\cK}_k(I- \mathbf{G}^k_{u}  \bm{\cK}_k)^{-1},\\
\bm{\nu}_k & := (I + \bm{\cL}_k \mathbf{G}^k_{u}) \bm{\upsilon}_k.
\end{align}
\end{subequations}

As is shown in \cite{p:bakCDC16,p:BAKOLAS2018}, the predicted cost
can be expressed as a convex function of the new decision variables
$(\bm{\cL}_k, \bm{\nu}_k )$; this new expression is denoted as
$\cJ(\bm{\cL}_k, \bm{\nu}_k)$. In addition, the constraint functions
$C_1(\bm{\cK}_k, \bm{\upsilon}_k)$ and $C_2(\bm{\cK}_k,
\bm{\upsilon}_k)$ become $\cC_1(\bm{\cL}_k, \bm{\nu}_k)$ and
$\cC_2(\bm{\cL}_k, \bm{\nu}_k)$, respectively. In particular,
$\cC_1(\bm{\cL}_k, \bm{\nu}_k)$ corresponds to an affine function in
$(\bm{\cL}_k, \bm{\nu}_k)$ whereas the constraint $\cC_2(\bm{\cL}_k,
\bm{\nu}_k) \in \mathbb{S}^{+}_n$ can be expressed as an LMI
constraint in terms of $(\bm{\cL}_k, \bm{\nu}_k)$ as is shown
in~\cite{p:bakCDC16,p:BAKOLAS2018,p:bakCDC2018}. The reader may
refer to the latter references for the technical details on the
conversion of the latter problems into tractable convex programs.

\subsection{Closed-Loop Nonlinear Dynamics and Propagation of Uncertainty}
Now let $\pi = \{\kappa(t,\cdot): t\in[0,N-1]\}$ be an admissible
control policy for Problem~\ref{problem1}. Then, the state space
model of the closed loop system is given by
\begin{equation}\label{eq:nonlcloop}
x(t+1) = f_{\mathrm{cl}}(t, x(t)) + w(t),
\end{equation}
where
\begin{equation}\label{eq:fcloop}
f_{\mathrm{cl}}(t,x) := f(x,\kappa(t,x)).
\end{equation}


Next, we describe the main steps for the propagation of the mean and
the covariance of the uncertain state of the nonlinear system
described by \eqref{eq:nonlcloop} based on the (scaled) unscented
transform~\cite{p:julier2002,p:julier2004}. To this aim, let us
assume that the mean $\mu_{k}:= \mathbb{E}[x(k)]$ and the covariance
$\Sigma_k:=\mathrm{Cov}[ x(k)]$ of the state of \eqref{eq:nonlcloop}
are known at stage $k$ (in practice only estimates / approximations
of the latter quantities will be known). Then, we will compute
$2n+1$ (deterministic) points, known as \textit{sigma points}, by
using the following equation:
\begin{align}\label{eq:sigma}
\sigma_k^{(i)} &= \begin{cases} \mu_k,&\mathrm{if}~i=0, \\
\mu_k + \sqrt{n + \lambda}
\Sigma_k^{1/2}\e_i, &\mathrm{if}~i\in [1,n]_{d}, \\
\mu_k - \sqrt{n + \lambda} \Sigma_k^{1/2}\e_{i-n},&\mathrm{if}~i\in
[n+1,2n]_{d},
\end{cases}
\end{align}
where $\{\e_i:~i \in [1,n]_{d}\}$ denotes the standard orthonormal
basis of $\mathbb{R}^{n}$. To each sigma point, we associate a pair
of gains $(\gamma_k^{(i)}, \delta_k^{(i)})$ where 
\begin{align}\label{eq:weight1}
\gamma_k^{(i)} & = \begin{cases} \lambda/(\lambda + n), &\mathrm{if}~i=0\\
1/(2( \lambda + n)), &\mathrm{if}~ i\in [1,2n]_{d},
\end{cases}
\end{align}
and
\begin{align}\label{eq:weight2}
\delta_k^{(i)} & = \begin{cases}
1 - \alpha^2 + \beta + \lambda/(\lambda + n), &\mathrm{if}~i=0,\\
1/(2( \lambda + n)), &\mathrm{if}~ i \in [1,2n]_{d}.
\end{cases}
\end{align}
The parameter $\alpha$ determines the spread around $\mu_k$ whereas
$\beta$ is a positive number and $\lambda:= \alpha^2 n -n $.
Typically, $0 < \alpha\ll 1$ and $\beta=2$ for Gaussian
approximations as suggested in \cite{p:wan2000,p:julier2004}.

Subsequently, we propagate the set of sigma points $\{
\sigma_k^{(i)}:~i \in [1, 2n+1]_{d} \}$ at the next stage $t=k+1$ to
obtain a new set of points $\{ \hat{\sigma}_{k+1}^{(i)}:~i \in [1,
2n+1]_{d} \}$, where
\begin{equation}\label{eq:sigma2}
\hat{\sigma}_{k+1}^{(i)} = f_{\mathrm{cl}}(k, \sigma_{k}^{(i)}),~~~i
\in [0,2n]_{d}.
\end{equation}
Using the point-set $\{ \hat{\sigma}_{k+1}^{(i)}:~i\in[0,2n]_{d}
\}$, one can approximate the (predicted) state mean and state
covariance at stage $t=k+1$ as follows:
\begin{subequations}
\begin{align}\label{eq:Xiupdate}
\hat{\mu}_{x}(k+1) & = \sum_{i=0}^{2L} \gamma_k^{(i)} \hat{\sigma}_{k+1}^{(i)},\\
\hat{\Sigma}_{x}(k+1) & = \sum_{i=0}^{2L} \delta_k^{(i)} (
\hat{\sigma}_{k+1}^{(i)} - \hat{\mu}_{x}(k+1) ) \nonumber \\
&~~~~~~\times ( \hat{\sigma}_{k+1}^{(i)} -
\hat{\mu}_{x}(k+1) )\t + W_k, \label{eq:Xiupdate2}
\end{align}
\end{subequations}
where $W_k\in \mathbb{S}^{+}_n$ corresponds to the noise covariance
at stage $t=k$.

\section{A greedy algorithm for nonlinear covariance
steering}\label{s:greedy}

The proposed algorithm consists of three main steps. We will
describe these steps starting at stage $t=k$, where $k\in[0,N-1]_d$,
and we will assume that approximations of the state mean
$\hat{\mu}_k$, the state covariance $\hat{\Sigma}_k$, and the input
mean $\hat{\nu}_k$ are known (if $k=0$, then we set $\hat{\mu}_k =
\mu_0$, $\hat{\Sigma}_k = \Sigma_0$, and $\hat{\nu}_k=0$).

We refer to the first step as the \textit{recursive linearization}
step (RL step). In the RL step, we construct a linearization $(A_k,
B_k, r_k)$ of \eqref{eq:nonlinear} around the point $(\hat{\mu}_k,
\hat{\nu}_k)$ by using \eqref{eq:AkBk}--\eqref{eq:rk}. Note that the
approximations $\hat{\mu}_k$ and $\hat{\nu}_k$ will be updated at
the end of each stage and consequently, the linearized model will
also have to be updated at each new stage to reflect the new
information and hence the ``recursive'' qualifier in the name of
this step. We will write
\begin{equation}
(A_k, B_k, r_k) = \Lambda \big[ \hat{\mu}_k, \hat{\nu}_k; f(\cdot)
\big].
\end{equation}

In the second step, which we refer to as the \textit{linearized
Gaussian covariance steering} step (LGCS step), we compute a
feedback control policy (sequence of feedback control laws) that
solves the $k$-th linearized covariance steering problem
(Problem~\ref{problem2}). To solve the latter problem, we need to
know the linearized model $(A_k,B_k,r_k)$, the approximations of the
predicted mean and covariance $(\hat{\mu}_k, \hat{\Sigma}_k)$ at
stage $k$ assuming that the goal state mean and state covariance
$(\mu_\f, \Sigma_\f)$ are known a priori. The triplet
$(A_k,B_k,r_k)$ is computed in the RL step, whereas the pair
$(\hat{\mu}_k, \hat{\Sigma}_k)$ is computed at the previous stage
(by executing the third step of the algorithm that will be discussed
shortly next). The policy $\varpi^{\star}_k$ that solves the $k$-th
linearized covariance steering problem 
with boundary conditions
\begin{subequations}
\begin{align}
\mathbb{E}[z(k)] & = \hat{\mu}_{k}, & \mathrm{Cov}[ z(k) ] =
\hat{\Sigma}_{k}, \\
\mathbb{E}[z(N)] & = \mu_\f, & (\Sigma_{\f} - \mathrm{Cov}[ z(N) ])
\in \mathbb{S}^{+}_n,
\end{align}
\end{subequations}
where $z$ corresponds to the state of the linearized system. We
write
\begin{equation}
\varpi_k^{\star} := \cS_k\big[ A_k, B_k, r_k, \hat{\mu}_{k},
\hat{\nu}_k, \hat{\Sigma}_{k} \big],
\end{equation}
where $\varpi_k^{\star} := \{ \phi^{\star}_k(k,\cdot), \dots,
\phi^{\star}_{k}(N-1,\cdot) \}$. The computation of the control
policy $\varpi_k^{\star}$ can be done in real-time by means of
robust and efficient convex optimization techniques (for details the
reader should refer to~\cite{p:bakCDC16,p:BAKOLAS2018}). We refer to
the latter step as the \textit{linearized Gaussian covariance
steering} step (LGCS step). After the computation of
$\varpi_k^{\star}$, we extract from it its first control law,
$\phi_k^{\star}(k,\cdot)$, that is, the control law that corresponds
to stage $k$. We write
\[
\phi_k^{\star}(k,z) := \cP_1\left[\varpi_k^{\star}\right] =
\upsilon^{\star}_k(k) + K^{\star}_k(k) z,
\]
where $\cP_1(\cdot)$ denotes the truncation operator that truncates
all the elements of a sequences except from the first one. Then, we
set the control law $\kappa^\star(k,\cdot)$ corresponding to the
$k$-th element of the feedback control policy $\pi^{\star}$ for the
original nonlinear covariance steering problem
(Problem~\ref{problem1}) to be equal $\phi_k^{\star}(k,\cdot)$, that
is,
\begin{equation}
\kappa^{\star}(k,x) := \phi_k^{\star}(k,x) = \upsilon^{\star}_k(k) +
K^{\star}_k(k) x,
\end{equation}
where $x$ is the state of the original nonlinear system.
Consequently, the one-stage transition map for the closed-loop
dynamics based on information available at stage $k$ is described by
the following equation:
\begin{equation}
x(k+1) = f^k_{\mathrm{cl}}(k,x) + w(k),
\end{equation}
where
\begin{align}\label{eq:fkcloop}
f^k_{\mathrm{cl}}(k,x) &:= f(x, \kappa^{\star}(k,x)) = f(x,
\upsilon^{\star}_k(k) + K^{\star}_k(k) x).
\end{align}

In the third step, we compute approximations $\hat{\mu}_x(k+1)$ and
$\hat{\Sigma}_x(k+1)$ of the (predicted) mean and covariance of the
state of the closed-loop system at stage $t=k+1$. To this aim, we
first compute a set of sigma points $\{
\hat{\sigma}_{k}:~k\in[0,2n]_d \}$ and their corresponding weights
$(\gamma^{(i)}_k, \delta^{(i)}_k)$ based on equations
\eqref{eq:sigma} and~\eqref{eq:weight1}-\eqref{eq:weight1},
respectively. Next, we compute the point-set $\{
\hat{\sigma}_{k+1}:~k\in[0,2n]_d \}$ by using
equation~\eqref{eq:sigma2} and the closed loop one-stage transition
map $f^k_{\mathrm{cl}}(k,x)$ which is defined in \eqref{eq:fkcloop}.
Subsequently, we compute $\hat{\mu}_{x}(k+1)$ and
$\hat{\Sigma}_{x}(k+1)$ by using
\eqref{eq:Xiupdate}-\eqref{eq:Xiupdate2}. The pair
$(\hat{\mu}_{x}(k+1),\hat{\Sigma}_{x}(k+1))$ determines a Gaussian
approximation of the statistics of the state of the closed-loop
system at stage $t=k+1$. We set $\hat{\mu}_{k+1} :=
\hat{\mu}_x(k+1)$ and $\hat{\Sigma}_{k+1} := \hat{\Sigma}_x(k+1)$.
Finally, we set $\hat{\nu}_{k+1} :=  \phi_k^{\star}(k+1,
\hat{\mu}_{k+1})$. We refer to the third step as the
\textit{predictive normalization} step (PN step). We write
\begin{equation}
(\hat{\mu}_{k+1}, \hat{\nu}_{k+1}, \hat{\Sigma}_{k+1}) :=
\mathcal{F}_k\big[ \hat{\mu}_{k},\hat{\Sigma}_{k};
f^k_{\mathrm{cl}}(k, \cdot) \big].
\end{equation}


These three steps of the previously described iterative process are
repeated for all stages $t \in [k,N-1]_d$ for a given $k\in
[0,N-1]_d$. At the end of the process, the predicted approximations
of the state mean and covariance are sufficiently close to their
corresponding goal quantities. The output of this iterative process
will be a control policy $\pi^{\star}_{k:N-1} :=\{
\kappa^{\star}(t,x):~t\in[k,N-1]_{\mZ}\}$. If $k \in [1,N-1]_{d}$,
then the policy $\pi^{\star}_{k:N-1}$ corresponds to the truncation
of the control policy $\pi^{\star}$ that solves
Problem~\ref{problem1}, which is comprised of the ``last'' $N-k$
elements of the latter policy. If we start the iterative process at
$k=0$, then the output of the process is the control policy
$\pi^{\star}$ that solves Problem~\ref{problem1}. The pseudocode of
the previous process is given in Algorithm~1.

\begin{algorithm}
\caption{Computation of feeback policy $\pi^{\star}_{k:N-1} :=
\{\kappa^{\star}(t,\cdot):t\in [k, N-1]_{\mZ}\}$ that solves
Problem~\ref{problem1}}\label{alg:iterncov}
\begin{algorithmic}[1]
\Procedure{Greedy Nonlinear Covariance Steering}{} 
\BState \textit{Input data}: $N$, $\mu_\f$, $\Sigma_\f$,
$f(\cdot)$ 
\BState \textit{Input variables}: $k$, $\hat{\mu}_{k}$, $\hat{\nu}_k$,
$\hat{\Sigma}_{k}$ 
\BState \textit{Output variables}: $\pi_{k:N-1}$,
$\{\hat{\mu}_x(t)\}_{t=k}^N$, $\{ \hat{\Sigma}_x(t) \}_{t=k}^N$
\For{$t=k:N-1$} \State $(A_t, B_t, r_t) := \Lambda\big[ \hat{\mu}_t,
\hat{\nu}_t, \hat{\Sigma}_{t} \big]$ \State $\varpi_t^{\star} :=
\cS_t[A_t, B_t, r_t, \hat{\mu}_{t}, \hat{\Sigma}_{t} ]$
\State $\kappa^{\star}(t,\cdot) := \cP_1[\varpi_t^{\star}]$ \State
$f^t_{\mathrm{cl}}(t,\cdot) := f(x, \kappa^{\star}(t,\cdot))$
\State$ (\hat{\mu}_{t+1}, \hat{\nu}_{t+1}, \hat{\Sigma}_{t+1}) :=
\mathcal{F}_t\big[
\hat{\mu}_{t},\hat{\Sigma}_{t};f^t_{\mathrm{cl}}(t,\cdot)  \big]$
\EndFor \State $\pi^{\star}_{k:N-1} := \{\kappa^{\star}(t,x):t\in
[k, N-1]_{\mZ}\}$ \EndProcedure
\end{algorithmic}
\end{algorithm}

\section{Numerical Simulations}\label{s:simu}

In this section, we present numerical simulations to illustrate the
basic ideas of this paper. In particular, we consider the following
DTSN system:
\begin{subequations}
\begin{align}
x_1(t+1) & = x_1(t) + \tau x_2(t),\\
x_2(t+1) & = x_2(t) - \tau ( \delta x_1(t) + \zeta x_1(t)^3 + \gamma x_2(t)  ) \nonumber \\
&~\quad~ + \tau u(t) + \sqrt{\tau} w(t),
\end{align}
\end{subequations}
where $[x_1(0),~x_2(0)]\t \sim \cN(\mu_0, \Sigma_0)$ with $\mu_0 =
[0,~0]\t$ and $\Sigma_0 = \mathrm{diag}(\sigma^2_1,\sigma_2^2)$,
where $\sigma_1=2.5$ and $\sigma_2 = 2.0$. In addition, the desired
terminal state mean and covariance are taken to be, respectively,
$\mu_\f = [0,~0]\t$ and $\Sigma_\f = \mathrm{diag}(s^2_1,s_2^2)$,
where $s_1 = 1.25$ and $s_2 = 1.0$. For our simulations, we consider the following parameter
values: $\tau = 0.01$, $N=100$ and $\zeta= 0.05$, $\gamma=0.05$, $\delta=-1$, $\alpha=0.05$, and $\beta=2$.

\begin{figure}
{ \epsfig{file =
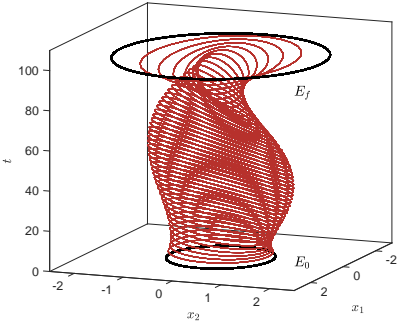,clip=,width=0.68\linewidth} \caption{\small{Time
evolution of the sequence $\{ E_t\}_{t=0}^{N}$. The vertical axis
in this 3D graph corresponds to the (discrete)
time-axis.}}\label{F:3Dplot}}
\end{figure}

Figure \ref{F:3Dplot} illustrates the time
evolution of the predicted state covariance $\hat{\Sigma}_x(t)$ in
terms of the evolution of the sequence of ellipsoids $\{
E_t\}_{t=0}^{N}$, where
\begin{equation*}
E_t : = \{ \mathsf{x} \in \mathbb{R}^2: (\mathsf{x}-
\hat{\mu}_x(t))\t \hat{\Sigma}_x(t)^{-1} (\mathsf{x} -
\hat{\mu}_x(t)) = 1 \},
\end{equation*}
for $t \in [0,N]_d$ (the ellipsoid $E_t$ is in an one-to-one
correspondence with $\hat{\Sigma}_x(t)$). To the desired terminal
state covariance $\Sigma_\f$, we associate the ellipsoid $E_\f$,
where
\begin{align*}
E_\f & : = \{ \mathsf{x} \in \mathbb{R}^2: ( \mathsf{x} - \mu_\f)\t
\Sigma_\f^{-1} ( \mathsf{x} -\mu_\f) =1  \}.
\end{align*}

\begin{figure}
\centering {
\epsfig{file = 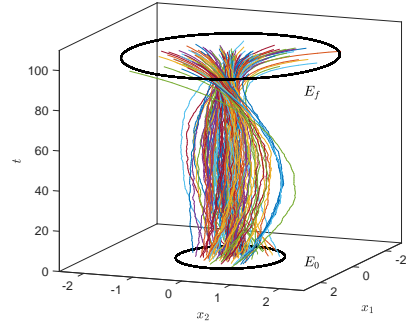,clip=,width=0.68\linewidth}
\caption{\small{Sample trajectories of the closed loop system.}}\label{F:trajectories}}
\end{figure}

\begin{figure}
\centering {
\epsfig{file = 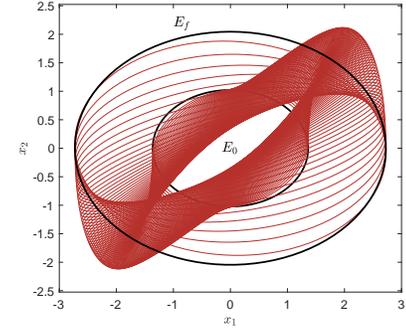,clip=,width=0.68\linewidth}
\caption{\small{Time evolution of $\{E_t\}_{t=0}^N$. The black ellipses correspond to $E_0$ and
$E_\f$.}}\label{F:2Dplot}}
\end{figure}

In particular, Fig.~\ref{F:3Dplot} illustrates the evolution of the
sequence $\{E_t\}_{t=0}^N$ in a 3D graph whose vertical axis
corresponds to the time-axis. Sample trajectories of the closed loop system are illustrated in Fig.~\ref{F:trajectories}. The projection on the $x_1-x_2$ plane of the 3D graph given in Fig.~\ref{F:3Dplot} is illustrated in Fig.~\ref{F:2Dplot}. In these three figures, the black ellipses correspond to $E_0$
and $E_\f$. We observe
that $E_N$ is very close to $E_\f$ and thus, the predicted
covariance of the terminal state $\hat{\Sigma}_\f$ is very close to
the goal covariance $\Sigma_\f$.

\begin{figure}
\centering  {
\epsfig{file = 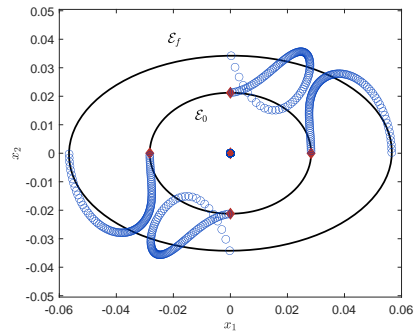,clip=,width=0.68\linewidth}
\caption{\small{Time evolution of the sigma points for
$t\in[0,N]_d$. The black ellipses correspond to the
ellipses $\cE_0$ and $\cE_\f$.}}\label{F:SigmaPoints}}
\end{figure}

The evolution of the sigma points used in the unscented transform
for the prediction of the state mean and covariance are illustrated
in Fig.~\ref{F:SigmaPoints}. In particular, the red diamonds
correspond to the sigma points associated with the initial state
mean and covariance, whereas the magenta circles correspond to the
predicted sigma points generated for all subsequent stages
$t\in[1,N]_d$. The sigma points corresponding to the original pair
$(\mu_0,\Sigma_0)$ and the terminal pair $(\mu_\f,\Sigma_\f)$ belong
to the ellipses $\cE_0$ and $\cE_\f$, respectively, where 
\begin{align*}
\cE_0 & : = \{ \mathsf{x} \in \mathbb{R}^2: (\mathsf{x}-\mu_0)\t
\Sigma_0^{-1} (\mathsf{x}-\mu_0)
= 2+\lambda \},\\
\cE_\f & : = \{ \mathsf{x} \in \mathbb{R}^2: (\mathsf{x}-\mu_\f)\t
\Sigma_\f^{-1} (\mathsf{x}-\mu_\f) = 2+\lambda  \}.
\end{align*}

\section{Conclusion}\label{s:concl}

In this work, we have proposed a greedy
covariance steering algorithm for discrete-time stochastic nonlinear
systems. The proposed approach relies on the solution of a sequence
of linearized covariance steering problems combined with the
(scaled) unscented transform that provides the one-stage predictions
of the mean and covariance of the state of the closed loop system.


To put the presented work under the umbrella of stochastic model
predictive control, it is necessary that performance and stability
considerations as well as notions of invariance based on
reachability analysis are integrated in the proposed algorithm. It
is worth noting that the reachability analysis for nonlinear
covariance steering problems requires the characterization of
``admissible'' sets of positive-definite matrices from which the
system can be steered to the desired state terminal covariance in
the given time horizon. To the best of our knowledge, the latter
reachability problem constitutes, at least for the case of
stochastic nonlinear systems, an open problem. In our future work,
we plan to study the latter problem and we will also explore
possible connections of this work with modern techniques of
stochastic model predictive control.

Another important problem in the context of nonlinear covariance
steering is the problem of verification of the results obtained with
the proposed greedy algorithm. At present, one can expect that the
predicted state mean and state covariance of the SNDT system, which
are computed by means of the unscented transform, will end up
sufficiently close to their goal quantities but this is not
automatically the case for the true state mean and state covariance
of the SNDT system. Finally, we plan to consider the case of
incomplete state information and also explore connections with
recent results on PDE tracking for distribution steering problems.

\bibliographystyle{ieeetr}
\bibliography{bibstocha}
\end{document}